\documentclass[a4paper, twoside]{scrartcl}
\usepackage[utf8]{inputenc}
\usepackage{scrlayer-scrpage, xhfill}
\usepackage{amsfonts,amssymb,amsmath,amsthm,amstext,amssymb,mathtools,mathrsfs,yfonts}
\usepackage{color,xcolor,enumitem,esint}
\usepackage{tikz}
\usepackage{dsfont}  % for indicator-one
\usepackage[authoryear,sort,round]{natbib}  % for bibliography
\usepackage{hyperref}
\hypersetup{
	breaklinks=true,
	colorlinks=true,
	linkcolor=blue,
	citecolor=blue,
	urlcolor=blue,
}

\usepackage[noabbrev,capitalise,nosort,nameinlink]{cleveref}

\crefname{assumption}{Assumption}{Assumptions}
\crefformat{footnote}{#2\footnotemark[#1]#3}

\addtokomafont{disposition}{\sffamily}

\usepackage{xpatch} % For QED symbol

\newcommand{\proofheadfont}{\bfseries\sffamily}
\xpatchcmd{\proof}{\itshape}{\proofheadfont}{}{}

\usepackage[hang,flushmargin]{footmisc}

\usepackage{mleftright}
\mleftright % Redefines \left to \mleft and \right to \mright globally

\DeclareMathOperator{\Span}{span}

\DeclareMathOperator{\argmin}{arg\,min}

\newcommand*{\defeq}{\coloneqq}

\renewcommand*{\emptyset}{\varnothing}
\renewcommand*{\epsilon}{\varepsilon}
\newcommand*{\cH}{\mathcal{H}}
\newcommand*{\bE}{\mathbb{E}}
\newcommand*{\bN}{\mathbb{N}}
\newcommand*{\bP}{\mathbb{P}}
\newcommand*{\cN}{\mathcal{N}}
\newcommand*{\cS}{\mathcal{S}}

\newcommand*{\cV}{\mathcal{V}}

\newcommand*{\qefed}{\eqqcolon}
\newcommand*{\quark}{\setbox0\hbox{$x$}\hbox to\wd0{\hss$\cdot$\hss}}
\newcommand*{\rd}{\mathrm{d}}

\newcommand*{\bR}{\mathbb{R}}

\DeclareMathOperator*{\supp}{supp}

\newcommand*{\perpspace}[1]{{#1}^{\perp}}
\newcommand*{\perpspaceball}[3]{B_{\perpspace{#1}}(#2, #3)}
\newcommand*{\oball}[2]{B_{\cH}(#1, #2)}
\newcommand*{\Oball}[2]{B_{\cH}\left(#1, #2\right)}

\newcommand*{\angball}[2]{B_{\cS}(#1, #2)}
\newcommand*{\Angball}[2]{B_{\cS}\left(#1, #2\right)}
\newcommand*{\cone}[2]{K_{\cH}(#1, #2)}
\newcommand*{\CMspace}[1]{H(#1)}

\renewcommand*{\geq}{\geqslant}

\renewcommand*{\leq}{\leqslant}

\renewcommand*{\mathbf}{\boldsymbol}

\newcounter{ntodo}
\newcommand{\todo}[1]{\stepcounter{ntodo}\bgroup\color{red}#1\egroup}
\newcommand{\reporttodo}{\bgroup\color{red}There are \arabic{ntodo} ``todo'' commands in the text.\egroup}

\newcounter{nHL}
\newcommand{\HLsays}[1]{\stepcounter{nHL}\bgroup\color{orange}Hefin:~#1\egroup}
\newcommand{\reportHL}{\bgroup\color{orange}There are \arabic{nHL} ``HL/Hefin'' commands in the text.\egroup}
\newcounter{nTJS}
\newcommand{\TJSsays}[1]{\stepcounter{nTJS}\bgroup\color{teal}Tim:~#1\egroup}
\newcommand{\reportTJS}{\bgroup\color{teal}There are \arabic{nTJS} ``TJS/Tim'' commands in the text.\egroup}
\newcounter{nST}
\newcommand{\STsays}[1]{\stepcounter{nST}\bgroup\color{brown}Sam:~#1\egroup}
\newcommand{\reportST}{\bgroup\color{brown}There are \arabic{nST} ``ST/Sam'' commands in the text.\egroup}

\theoremstyle{plain}
\newtheorem{theorem}{\sffamily Theorem}[section]
\newtheorem{proposition}[theorem]{\sffamily Proposition}
\newtheorem{lemma}[theorem]{\sffamily Lemma}
\newtheorem{corollary}[theorem]{\sffamily Corollary}
\theoremstyle{definition}
\newtheorem{definition}[theorem]{\sffamily Definition}

\newtheorem{example}[theorem]{\sffamily Example}

\newtheorem{remark}[theorem]{\sffamily Remark}

\newcommand{\absval}[1]{\lvert #1 \rvert}

\newcommand{\innerprod}[2]{\langle #1 , #2 \rangle}
\newcommand{\norm}[1]{\lVert #1 \rVert}
\newcommand{\set}[2]{\{ #1 \mid #2 \}}

\newcommand{\biginnerprod}[2]{\bigl\langle #1 , #2 \bigr\rangle}
\newcommand{\bignorm}[1]{\bigl\Vert #1 \bigr\Vert}
\newcommand{\bigset}[2]{\bigl\{ #1 \,\big\vert\, #2 \bigr\}}

\newcommand{\Set}[2]{\left\{ #1 \,\middle\vert\, #2 \right\}}

\numberwithin{equation}{section}
\numberwithin{figure}{section}
\numberwithin{table}{section}

\definecolor{darkgreen}{rgb}{0,0.6,0}

\usepackage{acro}
\DeclareAcronym{ACG}{short=ACG, long=angular central Gaussian}
\DeclareAcronym{BIP}{short=BIP, long=Bayesian inverse problem}
\DeclareAcronym{MAP}{short=MAP, long=maximum a posteriori}
\DeclareAcronym{OM}{short=OM, long=Onsager--Machlup}
\DeclareAcronym{SPSD}{short=SPSD, long=self-adjoint and positive semi-definite}
\definecolor{accentcol}{HTML}{E69F00}
\definecolor{accentcol2}{HTML}{56B4E9}
\definecolor{accentcol3}{HTML}{CC79A7}

\usepackage[a4paper, top=88pt, bottom=88pt, left=72pt, right=72pt, headsep=16pt, footskip=28pt]{geometry}

\newcommand*{\arXiv}[1]{\bgroup\color{blue}\href{https://arxiv.org/abs/#1}{arXiv:#1}\egroup}
\newcommand*{\doi}[1]{\bgroup\color{blue}\href{https://doi.org/#1}{doi:\nolinkurl{#1}}\egroup}
\newcommand*{\email}[1]{\bgroup\color{blue}\href{mailto:#1}{\nolinkurl{#1}}\egroup}
\renewcommand*{\url}[1]{\bgroup\color{blue}\href{#1}{\nolinkurl{#1}}\egroup}
\setlist[enumerate]{nosep}
\setlist[itemize]{nosep}

\automark[section]{section}
\setkomafont{pageheadfoot}{\normalcolor\sffamily}
\setkomafont{pagenumber}{\normalfont\normalsize\sffamily}
\clearpairofpagestyles
\let\oldtitle\title
\renewcommand{\title}[1]{\oldtitle{#1}\newcommand{\theshorttitle}{#1}}
\newcommand{\shorttitle}[1]{\renewcommand{\theshorttitle}{#1}}
\let\oldauthor\author
\renewcommand{\author}[1]{\oldauthor{#1}\newcommand{\theshortauthor}{#1}}
\newcommand{\shortauthor}[1]{\renewcommand{\theshortauthor}{#1}}
\cohead{\xrfill[0.525ex]{0.6pt}~\theshorttitle~\xrfill[0.525ex]{0.6pt}}
\cehead{\xrfill[0.525ex]{0.6pt}~\theshortauthor~\xrfill[0.525ex]{0.6pt}}
\newcommand{\theabstract}[1]{\par\bgroup\noindent\textbf{\textsf{Abstract.}} #1\egroup}
\newcommand{\thekeywords}[1]{\par\smallskip\bgroup\noindent\textbf{\textsf{Keywords.}}\newcommand{\and}{ $\bullet$ } #1\egroup}
\newcommand{\themsc}[1]{\par\smallskip\bgroup\noindent\textbf{\textsf{2020 Mathematics Subject Classification.}}\newcommand{\and}{ $\bullet$ } #1\egroup}
\newcommand*{\affilref}[1]{\ref{affiliation#1}}
\newcommand*{\affiliation}[3]{
	\footnotetext[#1]{\label{affiliation#2}#3}
}

\title{On the Onsager--Machlup functional of the infinite-dimensional angular central Gaussian distribution and related measures}
\shorttitle{The Onsager--Machlup functional of the angular central Gaussian distribution}
\author{%
	Hefin~Lambley\textsuperscript{\affilref{FUB},\affilref{WarwickSingle}}%
	\and%
	T.~J.~Sullivan\textsuperscript{\affilref{WarwickJoint}}%
	\and%
  Sam~Taylor\textsuperscript{\affilref{WarwickSingle}}%
}
\shortauthor{H.~Lambley, T.~J.~Sullivan, and S.~Taylor}
\date{\today}

\begin{document}

\maketitle
\affiliation{1}{FUB}{Institute of Mathematics, Freie Universit\"at Berlin, Arnimallee 6, 14195 Berlin, Germany\newline (\email{h.lambley@fu-berlin.de})}
\affiliation{2}{WarwickSingle}{Mathematics Institute, University of Warwick, Coventry, CV4 7AL, United Kingdom\newline (\email{hefin.lambley@warwick.ac.uk}, \email{samueltaylor1010@gmail.com})}
\affiliation{3}{WarwickJoint}{Mathematics Institute and School of Engineering, University of Warwick, Coventry, CV4 7AL, United Kingdom\newline (\email{t.j.sullivan@warwick.ac.uk})}

\begin{abstract}\small
	\theabstract{The Onsager--Machlup functional of a probability measure is a generalisation of the negative logarithm of the probability density function, defined precisely using the ratios of masses of balls in the limit of vanishing radius.
It plays an important role in several areas of applied mathematics because its minimisers can be understood as modes of the measure.
We show that the Onsager--Machlup functional of the angular central Gaussian measure on the unit sphere $\cS$ of an infinite-dimensional Hilbert space---i.e., the radial projection of a centred Gaussian measure with covariance operator $C$---is constant on the finite-dimensional subsphere of $\cS$ spanned by the eigenvectors of leading eigenvalue under $C$, and takes the value $+\infty$ elsewhere.
We discuss the consequences of this extreme concentration of mass for maximum a posteriori estimation in Bayesian inverse problems with an angular central Gaussian prior.
}
	\thekeywords{Angular central Gaussian measure%
\and%
Bayesian inverse problems%
\and%
maximum a posteriori estimation%
\and%
Onsager--Machlup functional%
%\and%
%\todo{More?}
}
	\themsc{28C20% Measure and integration > Set functions and measures on spaces with additional structure > Set functions and measures and integrals in infinite-dimensional spaces (Wiener measure, Gaussian measure, etc.)
\and%
58C35% Global analysis, analysis on manifolds > Calculus on manifolds; nonlinear operators > Integration on manifolds; measures on manifolds
\and%
60B05% Probability theory and stochastic processes > Probability theory on algebraic and topological structures > Probability measures on topological spaces
%\and%
%\todo{More?}
}
\end{abstract}

\section{Introduction}

The \ac{OM} functional of a probability measure acts as a kind of formal negative log-density for the measure, one that makes sense even in the absence of a uniform reference measure such as Lebesgue or Hausdorff measure.
It plays an important role in \ac{MAP} estimation, since its minimisers can be understood as modes, or most likely points, of the measure \citep{DuerrBach1978,AyanbayevKlebanovLieSullivan2022I}.

However, even some apparently simple measures lack well-understood \ac{OM} functionals.
This article offers a rigorous derivation of the \ac{OM} functional for the \ac{ACG} distribution---the radial projection of a centred Gaussian measure $\gamma_{0, C}$ on a Hilbert space $\cH$ onto its unit sphere $\cS$---as well as measures that are absolutely continuous with respect to the \ac{ACG} distribution.
The main result of this article, \Cref{thm:main}, shows that the \ac{OM} functional of the infinite-dimensional \ac{ACG} distribution is in fact highly degenerate:
it can be assigned any finite constant value on a particular finite-dimensional sub-sphere of $\cS$, and takes the value $+\infty$ elsewhere;
the sub-sphere in question is the intersection of $\cS$ with the span of the eigenvectors of the covariance operator $C$ that have maximal eigenvalue.

\Cref{sec:notation} sets out the notation and background concepts for the main result, \Cref{thm:main}, which is stated and proved in \Cref{sec:main}.
As well as proving \Cref{thm:main} in full rigour, we also provide the corresponding \ac{OM} functional in finite dimensions, which gives some intuition for the degenerate behaviour of the \ac{OM} functional in infinite dimensions.
Some comments on the consequences of \Cref{thm:main} for \ac{MAP} estimation in \acp{BIP} are given in \Cref{sec:discussion}, and some possible directions for further work are discussed in \Cref{sec:future}.

\section{Notation and background}
\label{sec:notation}

\subsection{General notation}

Given a separable metric space $(X, d)$, we denote the open ball in $X$ with centre $x \in X$ and radius $r > 0$ by $B(x, r) \defeq \set{y \in X}{d(x, y) < r}$.
Given a Borel probability measure $\mu$ on a separable metric space $X$, the \emph{topological support} $\supp(\mu)$ is
\begin{equation*}
	\supp(\mu) \defeq \bigset{ x \in X }{ \text{for all $r > 0$, } \mu\bigl(B(x, r)\bigr) > 0 }.
\end{equation*}
Throughout the paper, $\cH$ will be a separable real Hilbert space. %
The unit sphere $\cS$ of $\cH$ will be equipped with the \emph{angular metric} or \emph{geodesic metric} $d_{\cS}(u, v) \defeq \arccos \innerprod{ u }{ v }$.
Open balls in $\cH$ and $\cS$ will be denoted by
\begin{align*}
	\oball{v}{r} & \defeq \set{ u \in \cH }{ \norm{ u - v } < r } && \text{for $v \in \cH$, $r \geq 0$,} \\
	\angball{v}{\theta} & \defeq \set{ u \in \cS }{ d_{\cS}(u, v) < \theta } && \text{for $v \in \cS$, $\theta \geq 0$.}
\end{align*}
Note that a $d_{\cS}$-ball is a spherical cap, simultaneously the intersection of a half-space of $\cH$ with $\cS$ and the intersection of a ball in $\cH$ with $\cS$.
In particular,
\begin{equation*}
	\oball{v}{r} \cap \cS = \angball{v}{\theta}
	\quad
	\text{for $0 \leq r = 2 \sin \frac{\theta}{2} < 2$.}
\end{equation*}
We write $\cone{v}{\theta}$ for the open \emph{cone} with central ray $v \in \cH \setminus \{ 0 \}$ and aperture angle $2 \theta \geq 0$:
\[
	\cone{v}{\theta}
	\defeq \Set{ x \in \cH \setminus \{0\} }{ \arccos \left(\frac{ \innerprod{ x }{ v } }{ \norm{ x } \norm{ v } }\right) < \theta }
	= \bigcup_{\alpha > 0} \alpha \Angball{\frac{v}{\norm{ v }}}{\theta} .
\]
Given a separable Hilbert space $\cH$ and a closed subspace $\cV \subseteq \cH$, we denote by $\perpspace{\cV}$ the orthogonal complement of $\cV$, and by $P_{\cV} \colon \cH \to \cV$ the orthogonal projection onto $\cV$.
The inclusion $\iota \colon \cV \to \cH$ may be left implicit;
e.g., for $v \in \cV$ and $h \in \cH$, we write $v + h \in \cH$.
The adjoint of an operator $T \colon \cH \to \cV$ is denoted $T^{\ast} \colon \cV \to \cH$.
For vectors $x, y \in \cH$, we define the rank-one linear operator $(x \otimes y) \colon \cH \to \cH$ by
\begin{equation*}
	(x \otimes y)(h) \defeq \innerprod{y}{h} x \qquad \text{for $h \in \cH$.}
\end{equation*}
For a \ac{SPSD} and compact operator $C \colon \cH \to \cH$, we denote by $(\lambda_{n}(C))_{n \in \bN}$ the eigenvalues of $C$ in descending order, allowing for the possibility of repeated eigenvalues.
By the min-max principle \citep[Problem~37, p.~490]{Brezis2011}, the $n$\textsuperscript{th} eigenvalue $\lambda_{n}(C)$ can be written as follows, where the minimum is taken over all vector subspaces $W_{n-1} \subseteq \cH$ of dimension $n - 1$:
\begin{equation}
	\label{eq:Courant--Fischer--Weyl_min-max}
	\lambda_{n}(C) = \min_{\substack{W_{n-1} \subseteq \cH\\\dim W_{n-1} = n-1}} \max_{x \in W_{n-1}^{\perp} \cap \cS} \innerprod{C x}{x} .
\end{equation}

\subsection{Gaussian measures}

Let $\cH$ be an infinite-dimensional separable Hilbert space, let $m \in \cH$, and let $C \colon \cH \to \cH$ be a trace-class \ac{SPSD} operator.
Let $(\psi_{n})_{n \in \bN} \subset \cH$ be any sequence of orthonormal eigenvectors associated with the sequence $(\lambda_{n}(C))_{n \in \bN}$ of eigenvalues.
Then the \emph{Gaussian measure} $\gamma_{m, C}$ on $\cH$ with mean $m$ and covariance operator $C$ is the distribution of the random variable
\begin{equation}
	\label{eq:Karhunen--Loeve}
	x \defeq m + \sum_{n \in \bN} \sqrt{\lambda_{n}(C)} \psi_{n} \xi_{n}, \qquad \xi_{n} \overset{\text{i.i.d.}}{\sim} \cN(0, 1) \text{ in $\bR$.}
\end{equation}
The expression \eqref{eq:Karhunen--Loeve} is a \emph{Karhunen--Lo\`eve expansion} of $\gamma_{m, C}$ \citep[Theorem~3.5.1]{Bogachev1998}.
We say that $\gamma_{m, C}$ is \emph{non-degenerate} if it assigns positive mass to every open subset of $\cH$, or equivalently if the eigenvalues are all strictly positive, i.e., $C$ is positive definite.

The Anderson inequality for Gaussian measures \citep[Theorem~2.8.10]{Bogachev1998} states that, for any Gaussian measure $\gamma_{0, C}$ on a separable Hilbert space $\cH$,
\begin{equation}
	\label{eq:Anderson}
	\gamma_{0, C}\bigl( \oball{x}{r} \bigr) \leq \gamma_{0, C}\bigl(\oball{0}{r}\bigr) \quad \text{for all $x \in \cH$ and $r > 0$.}
\end{equation}

For each Gaussian measure $\gamma_{0, C}$, we define the \emph{Cameron--Martin space} $\CMspace{\gamma_{0, C}} \defeq  C^{1/2}(\cH)$.
When $\gamma_{0, C}$ is non-degenerate, $\CMspace{\gamma_{0, C}}$ is a dense subspace of $\cH$;
we further equip $\CMspace{\gamma_{0, C}}$ with an inner product, and corresponding induced norm, given by the following reweightings of the ambient inner product and norm in $\cH$:
\begin{equation*}
	\innerprod{h}{k}_{\CMspace{\gamma_{0, C}}} \defeq \biginnerprod{C^{-1/2} h}{C^{-1/2} k} , \qquad \norm{h}_{\CMspace{\gamma_{0, C}}} \defeq \bignorm{C^{-1/2} h} .
\end{equation*}
The Cameron--Martin space influences several measure-theoretic properties of the Gaussian measure $\gamma_{0, C}$, as will be discussed further in \Cref{subsec:infinite-dimensional}.

\subsection{Onsager--Machlup functionals}

\begin{definition}
	Let $X$ be a separable\footnote{Separability ensures that the Borel $\sigma$-algebra coincides with the $\sigma$-algebra generated by the open balls of $X$, and also that, for every $r > 0$, at least one open ball $B(x, r)$, $x \in X$, has strictly positive mass under $\mu$.} metric space and let $\mu$ be a probability measure on the Borel $\sigma$-algebra of $X$.
	A function $I_{\mu} \colon \emptyset \neq E \subseteq \supp(\mu) \to \bR$ is called an \emph{\acl{OM} functional} for $\mu$ if, for every $x, y \in E$,
	\[
		\lim_{r \to 0} \frac{\mu(B(x, r))}{\mu(B(y, r))} = \exp \bigl( - I_{\mu}(x) + I_{\mu}(y) \bigr) .
	\]
	The functional $I_{\mu}$ is called \emph{exhaustive}\footnote{This term is used, for example, by \citet{Lambley2026} and \citet{KlebanovLambleySullivan2025}.
	The property is also known as property $M(\mu, E)$ \citep{AyanbayevKlebanovLieSullivan2022I}, but we do not use this term here.
	Unnamed, it plays an important role in earlier treatments of \ac{MAP} estimation in function spaces \citep{DashtiLawStuartVoss2013}.}  if, for every $x \notin E$ and $y \in E$,
	\[
		\lim_{r \to 0} \frac{\mu(B(x, r))}{\mu(B(y, r))} = 0.
	\]
	An exhaustive \ac{OM} functional $I_{\mu} \colon E \to \bR$ may be extended to have domain $X$ by assigning the value $I_{\mu}(x) = + \infty$ for $x \notin E$.
	For brevity, we also use the term \emph{exhaustive \ac{OM} functional} to refer to such an extension $I \colon X \to \bR \cup \{+\infty\}$.%
\end{definition}

Note that an \ac{OM} functional is not unique, since if $I_{\mu}$ is an \ac{OM} functional for $\mu$, then so too is $I_{\mu} + c$, for any constant $c \in \bR$.
It is also important to note that many measures do not admit globally defined \ac{OM} functionals with $E = X$.
The standard example of this is a Gaussian measure $\gamma_{0, C}$ on a Hilbert space $\cH$:

\begin{example}[{\Ac{OM} functional for Gaussian measures; \citealp[Corollary~4.7.8]{Bogachev1998}}]
	\label{ex:Gaussian_OM}
	Let $\gamma_{0, C}$ be a centred, non-degenerate Gaussian measure on the separable Hilbert space $\cH$.
	Then $\gamma_{0, C}$ has an exhaustive \ac{OM} functional $I_{\gamma_{0, C}} \colon \cH \to \bR \cup \{+\infty\}$ that is finite on the Cameron--Martin space $E = \CMspace{\gamma_{0, C}} \defeq C^{1/2} (\cH)$, defined by
	\[
		I_{\gamma_{0, C}}(x) \defeq \frac{1}{2} \norm{ x }_{\CMspace{\gamma_{0, C}}}^{2} \defeq \frac{1}{2} \norm{ C^{-1/2} x }^{2}.
	\]
\end{example}

\subsection{Angular central Gaussian measures}

This paper's central object of interest, the \ac{ACG} measure, is defined as follows:

\begin{definition}
	Let $\gamma_{0, C}$ be a centred Gaussian measure on a separable Hilbert space $\cH$, with covariance operator $C \neq 0$.
	The \emph{angular central Gaussian measure} $\nu_{C}$ on the unit sphere $\cS$ of $\cH$ is then defined to be the push-forward of $\gamma_{0, C}$ via the radial projection map, i.e.,
	\[
		\nu_{C}(E) \defeq \gamma_{0, C} \left( \bigcup_{\alpha > 0} \alpha E \right)
		\quad
		\text{for all Borel-measurable $E \subseteq \cS$.}
	\]
	In particular, $\nu_{C} ( \angball{v}{\theta} ) = \gamma_{0, C} ( \cone{v}{\theta} )$.
\end{definition}

In directional statistics, the \ac{ACG} distribution models antipodally symmetric directional data.
It is an alternative to the Bingham distribution \citep{Tyler1987} and a special case of the projected normal distribution, for which the initial Gaussian distribution may have nonzero mean \citep{WangGelfand2013}.
\Ac{ACG} distributions and their mixtures have been applied in finite-dimensional directional-statistical problems such as geomagnetism \citep[Section~8]{Tyler1987}, imaging in neuroscience \citep{Tabelow2012}, and materials science \citep[Section~4]{Franke2016};
infinite-dimensional applications include binary classification and distributional inference problems in which one wishes to learn a \emph{normalised} function on some set \citep{LieRudolfSprungkSullivan2023}.

We conclude this section with a straightforward characterisation of what it means for an \ac{ACG} measure to be \emph{non-degenerate}:

\begin{lemma}
	\label[lemma]{lem:positivity}
	Let $\gamma_{0, C}$ be a centred Gaussian measure on a separable Hilbert space $\cH$ with non-zero covariance operator $C$ and let $\nu_{C}$ be the induced \ac{ACG} measure on the unit sphere $\cS$.
	Then the following are equivalent:
	\begin{enumerate}[label=(\alph*)]
		\item
		\label{item:positivity_cov}
		$C$ is positive definite;
		\item
		\label{item:positivity_Gaussian}
		$\gamma_{0, C}$ assigns strictly positive mass to every open subset of $\cH$, i.e., $\supp(\gamma_{0, C}) = \cH$;
		\item
		\label{item:positivity_ACG}
		$\nu_{C}$ assigns strictly positive mass to every open subset of $\cS$, i.e., $\supp(\nu_{C}) = \cS$.
	\end{enumerate}
\end{lemma}

\begin{proof}
	(\ref{item:positivity_cov}$\iff$\ref{item:positivity_Gaussian})~
	This equivalence is classical, and follows from the fact that $\supp(\gamma_{0, C})$ is the closure of $\CMspace{\gamma_{0, C}} = C^{1/2}(\cH)$ in $\cH$;
	see, e.g., \citet[Theorem~3.6.1]{Bogachev1998}.

	\smallskip\noindent(\ref{item:positivity_Gaussian}$\implies$\ref{item:positivity_ACG})~
	Suppose that \ref{item:positivity_Gaussian} holds, and let $U \subset \cS$ be open with respect to $d_{\cS}$.
	Let $u \in U$ and let $0 < \theta < \frac{\pi}{2}$ be such that $\angball{u}{\theta} \subseteq U$.
	Then $\oball{u}{\sin \theta} \subset \cone{u}{\theta}$, and hence
	\[
		\nu_{C}(U) \geq \nu_{C}(\angball{u}{\theta}) = \gamma_{0, C}(\cone{u}{\theta}) \geq \gamma_{0, C}(\oball{u}{\sin \theta}) > 0 .
	\]

	\smallskip\noindent(\ref{item:positivity_ACG}$\implies$\ref{item:positivity_cov})~
	If \ref{item:positivity_cov} fails, then there is some $u \in \cS \cap \ker C$.
	Then, for $x \sim \gamma_{0, C}$, $\bE [ \absval{ \innerprod{ x }{ u } }^{2} ] = \innerprod{ C u }{ u } = 0$, and so $x \in \perpspace{\{ u \}}$ $\gamma_{0, C}$-almost surely.
	But then \ref{item:positivity_ACG} fails since, for $0 < \theta < \tfrac{\pi}{2}$,
	\[
		\nu_{C}(\angball{u}{\theta}) = \gamma_{0, C}(\cone{u}{\theta}) \leq \gamma_{0, C} \bigl( \cH \setminus \perpspace{\{ u \}} \bigr) = 0 . \qedhere
	\]
\end{proof}

\section{Main result}
\label{sec:main}

\begin{theorem}
	\label{thm:main}
	Let $\nu_{C}$ be a non-degenerate \ac{ACG} measure on the unit sphere $\cS$ of an infinite-dimensional, separable Hilbert space $\cH$.
	Let $L = \ker(C - \lambda_{1}(C) I)$ be the subspace of $\cH$ spanned by the eigenvectors of $C$ with maximal eigenvalue.
	Then $\nu_{C}$ has an exhaustive \ac{OM} functional $I_{\nu_{C}} \colon \cS \to \bR \cup \{ + \infty \}$ given by
	\begin{equation} \label{eq:main}
		I_{\nu_{C}}(x) =
		\begin{cases}
			0, & \text{if $x \in \cS \cap L$,} \\
			+\infty, & \text{if $x \in \cS \setminus L$.}
		\end{cases}
	\end{equation}
\end{theorem}

\begin{remark}
	By Sazonov's theorem, the covariance operator $C \colon \cH \to \cH$ is a trace-class \ac{SPSD} operator.
	In particular, the eigenvalues of $C$ form a nonnegative null sequence, and the assumption that the measure is non-degenerate implies that all eigenvalues are strictly positive.
	Hence, $L$ always has strictly positive and finite dimension.
\end{remark}

\subsection{Finite-dimensional calculations}

Before proving \Cref{thm:main}, it is instructive to consider its finite-dimensional analogue, which also describes the \ac{OM} functional of an \ac{ACG} measure when the initial Gaussian is degenerate and the underlying covariance operator has finite rank $d$.

\begin{lemma}
	\label{lem:main_but_Rn}
	Let $\nu_{C}$ be a non-degenerate \ac{ACG} measure on the unit sphere $\cS \subset \bR^{d}$.
	Then $\nu_{C}$ has an exhaustive \ac{OM} functional $I_{\nu_{C}} \colon \cS \to \bR$ given by
	\begin{equation}
		\label{eq:finite-dimensional_ACG_OM}
		I_{\nu_{C}}(x) = d \log \sqrt{ \innerprod{x}{C^{-1} x} } .
	\end{equation}
\end{lemma}

\begin{proof}
	As noted by \citet[Eq.~(1)]{Tyler1987}, $\nu_{C}$ has a probability density function $f$ with respect to the $(d - 1)$-dimensional Hausdorff measure $\chi^{d - 1}$ on the sphere $\cS$ given by
	\[
		f(u) = \frac{\Gamma(d/2)}{2 \sqrt{ \pi^{d} \det C }} \innerprod{u}{C^{-1} u}^{-d / 2}.
	\]
	Observe that $f$ is a uniformly continuous function on $\cS$;
	let $\varpi_{f}$ be a modulus of continuity for $f$.
	Then, for $u, v \in \cS$,
	\begin{align*}
		\limsup_{\theta \to 0} \frac{\nu_{C}\bigl(\angball{u}{\theta}\bigr)}{\nu_{C}\bigl(\angball{v}{\theta}\bigr)}
		& = \limsup_{\theta \to 0} \frac{\int_{\angball{u}{\theta}} f(z) \, \chi^{d - 1}(\rd z)}{\int_{\angball{v}{\theta}} f(z) \, \chi^{d - 1}(\rd z)} \\
		& \leq \limsup_{\theta \to 0} \frac{\bigl(f(u) + \varpi_{f}(\theta)\bigr) \chi^{d - 1}\bigl(\angball{u}{\theta}\bigr)}{\bigl(f(v) - \varpi_{f}(\theta)\bigr) \chi^{d - 1}\bigl(\angball{v}{\theta}\bigr)} %
		= \frac{f(u)}{f(v)} ,
	\end{align*}
	since $\chi^{d - 1}(\angball{u}{\theta}) = \chi^{d - 1}(\angball{v}{\theta})$ and $\lim_{\theta \to 0} \varpi_{f}(\theta) = 0$.
	A similar calculation shows that
	\[
		\liminf_{\theta \to 0} \frac{\nu_{C}\bigl(\angball{u}{\theta}\bigr)}{\nu_{C}\bigl(\angball{v}{\theta}\bigr)}
		\geq \frac{f(u)}{f(v)} ,
	\]
	and combining these two facts shows that $- \log f$, which differs from \eqref{eq:finite-dimensional_ACG_OM} only by an additive constant, is an \ac{OM} functional for $\nu_{C}$.
\end{proof}

\begin{remark}
	\begin{enumerate}[label=(\alph*)]
		\item
		In view of \Cref{lem:main_but_Rn}, one might suspect that, as we consider Hilbert spaces of dimension $d \to \infty$, the \ac{OM} functionals of the corresponding \ac{ACG} measures become increasingly extreme:
		they either take their minimum value, which occurs on $\cS \cap L$, and this minimum value may be taken to be zero;
		or else they take a positive value of order $d$.
		Informally taking the limit as $d \to \infty$, \Cref{thm:main} is thus a plausible conjecture for an \ac{OM} functional of an infinite-dimensional \ac{ACG} measure.
		However, each application of \Cref{lem:main_but_Rn} considers a different measure on a different space, and so to establish \Cref{thm:main} rigorously we must analyse the infinite-dimensional \ac{ACG} measure directly.

		\item
		\Cref{lem:main_but_Rn} is easily adapted to show that, whenever a measure $\mu$ has a continuous density $f$ with respect to a strictly positive, locally finite, and uniformly distributed measure $\Upsilon$ on a separable metric space $X$---i.e., a measure such that the $\Upsilon$-mass of a ball is determined by its radius---it follows that $- \log f$ is an \ac{OM} functional for $\mu$.
		The absence of such a measure $\Upsilon$ on the infinite-dimensional sphere is another reason why the infinite-dimensional \ac{ACG} measure must be analysed directly and at more length than the proof of \Cref{lem:main_but_Rn}.
	\end{enumerate}
\end{remark}

\subsection{Proof of \texorpdfstring{\Cref{thm:main}}{Theorem~\ref{thm:main}}}
\label{subsec:infinite-dimensional}

\begin{proof}[Proof of \Cref{thm:main}]
	Let $(\psi_{n})_{n \in \bN}$ be an orthonormal eigenbasis of $\cH$ associated with the decreasing sequence $(\lambda_{n}(C))_{n \in \bN}$ of eigenvalues of $C$.
	It suffices to show that
	\begin{equation}
		\lim_{\theta \to 0} \frac{\nu_{C}(\angball{v}{\theta})}{\nu_{C}(\angball{\psi_{1}}{\theta})} = \begin{cases}
			1, & \text{for all $v \in \cS\cap L$},\\
			0, & \text{for all $v \in \cS \setminus L$}.
		\end{cases}
		\label{eq:proof_sketch_statement}
	\end{equation}
	Since only the behaviour as $\theta \to 0$ is relevant, we assume in what follows that $0 < \theta < \frac{\pi}{2}$.
	Recall that, by definition, $\nu_{C}(\angball{v}{\theta})$ is precisely the Gaussian cone mass $\gamma_{0, C}(\cone{v}{\theta})$.

	\paragraph{Step 1: Disintegrate $\gamma_{0, C}$ along the central ray of the cone.}
	Let $x \sim \gamma_{0, C}$, and write
	\begin{equation*}
		x = w + \tau v, \quad \text{where $w \defeq P_{\perpspace{\{v\}}} x$ and $\tau \defeq \innerprod{x}{v}$.}
	\end{equation*}
	Since the map $\cH \ni x \mapsto (w, \tau v) \in \perpspace{\{v\}} \oplus \Span\{v\}$ is linear, both $\tau$ and $w \mid \tau$ are Gaussian, and the standard conditioning formula for jointly Gaussian random variables in Hilbert space \citep[see, e.g.,][Lemma~4.3]{HairerStuartVossWiberg2005} implies that $\tau \defeq \innerprod{x}{v} \sim \gamma_{0,  \sigma_{v}^{2} }$ in $\bR$, for $\sigma^{2}_{v} \defeq \innerprod{Cv}{v}$, and
	$w \mid \tau$ has Gaussian distribution $\gamma_{\tau m_{v}, C_{v}}$ in $\perpspace{\{v\}}$ with parameters
	\begin{align}
		m_{v} &\defeq \sigma_{v}^{-2} P_{\perpspace{\{v\}}}  Cv,
		\label{eq:disintegration_mean} \\
	C_{v} &\defeq P_{\perpspace{\{v\}}} R_{v}  P^{\ast}_{\perpspace{\{v\}}}, \label{eq:disintegration_covariance} \\
		R_{v} &\defeq C - \sigma_{v}^{-2} \bigl( Cv \otimes Cv \bigr) . \label{eq:disintegration_covariance_R}
	\end{align}
	We write $w \mid \tau$ in terms of a Karhunen--Lo\`eve expansion \eqref{eq:Karhunen--Loeve}, with $(\varphi_{n})_{n \in \bN}$ denoting an orthonormal basis of $\perpspace{\{v\}}$ associated with the decreasing sequence $(\lambda_{n}(C_{v}))_{n \in \bN}$ of eigenvalues:
	\begin{equation*}
		w \mid \tau = \tau m_{v} + \sum_{n \in \bN} \sqrt{\lambda_{n}(C_{v})} \xi_{n} \varphi_{n},\qquad \xi_{n} \overset{\text{i.i.d.}}{\sim} \cN(0, 1).
	\end{equation*}
	Then, denoting by $\rho_{m, s^{2}}$ the density of a Gaussian on $\bR$ with mean $m$ and variance $s^{2}$, and by $\perpspaceball{\{v\}}{0}{t \tan \theta}$ the ball in $\perpspace{\{v\}}$ centred at $0$ with radius $t \tan \theta$, we have
	\begin{align}
		\nu_{C}\bigl(\angball{v}{\theta} \bigr) &= \bP\bigl[ x \in \cone{v}{\theta} \bigr] \nonumber\\
		&= \int_{0}^{\infty} \bP\left[ w \in \perpspaceball{\{v\}}{0}{t \tan \theta} \middle| \tau = t \right] \rho_{0,\sigma_{v}^{2}}(t) \,\rd t \nonumber \\
		&= \int_{0}^{\infty} \bP\left[ \sum_{n \in \bN} \sqrt{\lambda_{n}(C_{v})} \xi_{n} \varphi_{n} \in \perpspaceball{\{v\}}{-tm_{v}}{t \tan \theta} \right] \rho_{0, \sigma_{v}^{2}}(t) \,\rd t.
		\label{eq:nu_mass_to_disintegration}
	\end{align}

	\paragraph{Step 2: For $v \in \cS \cap L$, show that $\nu_{C}( \angball{v}{\theta} ) = \nu_{C}( \angball{\psi_{1}}{\theta} )$.}
	If $v \in \cS \cap L$, it is easy to verify that $\sigma_{v} = \sigma_{\psi_{1}}$ and $m_{v} = 0$.
	Moreover, the operators $R_{v}$ and $R_{\psi_{1}}$ are isospectral, i.e., $\lambda_{n}(R_{v}) = \lambda_{n}(R_{\psi_{1}})$ for all $n \in \bN$:
	by the definition \eqref{eq:disintegration_covariance_R}, the common spectrum of both operators is given by removing one copy of the leading eigenvalue from the spectrum of $C$.
	Since $R_{v} = C_{v} \oplus 0$ on $\cH = \perpspace{\{v\}} \oplus \Span\{v\}$, it follows that $\lambda_{n}(C_{v}) = \lambda_{n}(R_{v}) = \lambda_{n}(R_{\psi_{1}}) = \lambda_{n}(C_{\psi_{1}})$ for all $n \in \bN$.
	Thus, using \eqref{eq:nu_mass_to_disintegration}, we find
	\begin{equation} \label{eq:mass_same_in_leading_eigenspace}
		\nu_{C}\bigl( \angball{v}{\theta} \bigr) = 
		\int_{0}^{\infty} \bP\left[ \sum_{n \in \bN} \lambda_{n}(C_{\psi_{1}}) \xi_{n}^{2} < (t \tan \theta)^{2} \right] 
		\rho_{0, \sigma_{\psi_{1}}^{2}}(t) \,\rd t = \nu_{C}\bigl( \angball{\psi_{1}}{\theta} \bigr).
	\end{equation}
	This yields the claim \eqref{eq:proof_sketch_statement} for $v \in \cS \cap L$.

	\paragraph{Step 3: For $v \in \cS \setminus L$, show that $\nu_{C}(\angball{v}{\theta})$ is negligible relative to $\nu_{C}(\angball{\psi_{1}}{\theta})$ as $\theta \to 0$.}
	First, let $v \in \cS$ and $\theta > 0$.
	Applying the Anderson inequality \eqref{eq:Anderson} to \eqref{eq:nu_mass_to_disintegration}, then applying the min-max principle to the operator $C_{v}$, which yields that $\lambda_{n}(C_{\psi_{1}}) \leq \lambda_{n}(C_{v})$ for all $n \in \bN$ (\Cref{lem:min-max}), we find that
	\begin{align}
		\nu_{C}\bigl( \angball{v}{\theta} \bigr) 
		&\leq 
		\int_{0}^{\infty} \bP\left[ \sum_{n \in \bN} \sqrt{\lambda_{n}(C_{v})} \xi_{n} \varphi_{n} \in \perpspaceball{\{v\}}{0}{t \tan \theta} \right] \rho_{0, \sigma_{v}^{2}}(t) \,\rd t \nonumber \\
		&\leq 
		\int_{0}^{\infty} \bP\left[ \sum_{n \in \bN} \lambda_{n}(C_{\psi_{1}}) \xi_{n}^{2} < (t \tan \theta)^{2} \right] 
		\rho_{0, \sigma_{v}^{2}}(t)
		\,\rd t \label{eq:acg_exhaustivity_bound_1}  \\
		&\leq \int_{0}^{\infty} \bP\left[\sum_{n \in \bN} \lambda_{n}(C_{\psi_{1}}) \xi_{n}^{2} < (t \tan \theta)^{2} \right]  \rho_{0, \sigma_{\psi_{1}}^{2}}(t) \,\rd t = \nu_{C}\bigl( \angball{\psi_{1}}{\theta} \bigr). \label{eq:acg_exhaustivity_bound_2}
	\end{align}
	We use this to prove the following stronger statement:
	for any $0 < \varepsilon < \sigma_{\psi_{1}}$, any positive null sequence $(\theta_{j})_{j \in \bN}$, and any sequence $(v_{j})_{j \in \bN} \subset \cS$ satisfying $\sigma_{v_{j}} < \sigma_{\psi_{1}} - \varepsilon$ for all $j \in \bN$, we have
	\begin{equation}
		\label{eq:acg_exhaustivity}
		\lim_{j \to \infty} \frac{\nu_{C}\bigl( \angball{v_{j}}{\theta_{j}} \bigr)}{\nu_{C}\bigl( \angball{\psi_{1}}{\theta_{j}} \bigr)} = 0.
	\end{equation}
	This suffices to prove \eqref{eq:proof_sketch_statement} because, for every positive null sequence $(\theta_{j})_{j \in \bN}$, we may take the constant sequence $v_{j} = v$, and for $v \in \cS \setminus L$ we have that $\sigma_{v} < \sigma_{\psi_{1}}$.
	To show \eqref{eq:acg_exhaustivity}, consider \eqref{eq:acg_exhaustivity_bound_1}, write $F(r) \defeq \bP[\sum_{n \in \bN} \lambda_{n}(C_{\psi_{1}}) \xi_{n}^{2} < r^{2}]$, 
	and observe that, for any $M > 0$, we have
	\begin{align}
		\nu_{C}\bigl( \angball{v_{j}}{\theta_{j}} \bigr) &\leq %
		\int_{0}^{M} F(t \tan \theta_{j})
		\rho_{0, \sigma_{v_{j}}^{2}}(t)
		\,\rd t  \nonumber\\
		& \phantom{=} \quad 
		+ 
		\frac{\sigma_{\psi_{1}}}{\sigma_{v_{j}}} \int_{M}^{\infty} F(t \tan \theta_{j}) 
		\rho_{0, \sigma_{\psi_{1}}^{2}}(t) 
		\exp \left(\frac{t^{2}}{2} \left( \sigma_{\psi_{1}}^{-2} - \sigma_{v_{j}}^{-2} \right) \right) \,\rd t \nonumber  \\
		&\leq F(M \tan \theta_{j})
		+ \frac{\sigma_{\psi_{1}}}{\sigma_{v_{j}}} \exp\left( - \frac{g_{j} M^{2}}{2} \right) \nu_{C}\bigl( \angball{\psi_{1}}{\theta_{j}} \bigr) , \label{eq:nu_upper_bound}
	\end{align}
	where $g_{j} \defeq \sigma_{v_{j}}^{-2} - \sigma_{\psi_{1}}^{-2} > 0$.
	Meanwhile, writing the Gaussian density $\rho_{0, \sigma_{\psi_{1}}^{2}}$ explicitly yields
	\begin{align}
		\nu_{C}\bigl(\angball{\psi_{1}}{\theta_{j}} \bigr) &\geq 
		\int_{2M}^{3M}  F(t \tan \theta_{j}) 
		\rho_{0, \sigma_{\psi_{1}}^{2}}(t)
		\,\rd t %
		\geq \frac{M}{\sqrt{2\pi} \sigma_{\psi_{1}}} F(2M \tan \theta_{j}) \exp\left( -\frac{9M^{2}}{2 \sigma_{\psi_{1}}^{2}} \right). \label{eq:nu_lower_bound}
	\end{align}
	Therefore, combining the upper bound \eqref{eq:nu_upper_bound} and the lower bound \eqref{eq:nu_lower_bound},
	\begin{equation} \label{eq:ball_ratio_two-term_bound}
		\frac{\nu_{C}(\angball{v_{j}}{\theta_{j}})}{\nu_{C}(\angball{\psi_{1}}{\theta_{j}})} \leq \frac{\sqrt{2\pi} \sigma_{\psi_{1}}}{M} \exp\left(\frac{9M^{2}}{2\sigma_{\psi_{1}}^{2}}\right) \frac{F(M \tan \theta_{j})}{F(2M \tan \theta_{j})} + \frac{\sigma_{\psi_{1}}}{\sigma_{v_{j}}} \exp\left( -\frac{g_{j} M^{2}}{2} \right).
	\end{equation}
	By \Cref{lem:Gaussian_rapid_variation}, the Gaussian ball mass $F(r)$ is rapidly varying when $\cH$ has infinite dimension, meaning that the first term on the right-hand side of \eqref{eq:ball_ratio_two-term_bound} converges to zero as $j \to \infty$, i.e., as radius $\theta_{j} \to 0$.
	For the second term in \eqref{eq:ball_ratio_two-term_bound}, observe that, by hypothesis, $g_{j} \geq (\sigma_{\psi_{1}} - \varepsilon)^{-2} - \sigma_{\psi_{1}}^{-2} \qefed g^{\star} > 0$ for all $j \in \bN$.
	Thus, using the elementary bound $\log x \leq x - 1$ for $x > 0$, we find
	\begin{align*}
		\frac{\sigma_{\psi_{1}}}{\sigma_{v_{j}}} \exp\left( - \frac{g_{j} M^{2}}{2} \right)
		\leq \exp\left(\frac{1}{2} \left( \frac{\sigma_{\psi_{1}}^{2}}{\sigma_{v_{j}}^{2}} - 1 - g_{j} M^{2} \right) \right)
		\leq \exp\left( - \frac{(M^{2} - \sigma_{\psi_{1}}^{2}) g^{\star}}{2} \right).
	\end{align*}
	Since $g^{\star} > 0$ and the above bound holds for arbitrarily large $M > 0$, the right-hand side can be made arbitrarily small.
	This establishes \eqref{eq:acg_exhaustivity}, and the claim \eqref{eq:proof_sketch_statement} for $v \in \cS \setminus L$ follows.
\end{proof}

The proof of \Cref{thm:main} requires control of the Gaussian ball mass $\gamma_{0, C_{\psi_{1}}}(\oball{0}{r})$, written in terms of the Karhunen--Lo\`eve expansion as $F(r) \defeq \bP\bigl[ \sum_{n \in \bN} \lambda_{n}(C_{\psi_{1}}) \xi_{n}^{2} < r^{2} \bigr]$, in the limit $r \to 0$.
Before stating and proving this asymptotic, we recall an important property of the Cameron--Martin space---that it is the \emph{shift quasi-invariance space} of $\gamma_{0, C}$.
This means that, for any $h \in \CMspace{\gamma_{0, C}}$, the shifted Gaussian $\gamma_{h, C}$ is absolutely continuous with respect to $\gamma_{0, C}$.
Moreover, the density of $\gamma_{h, C}$ with respect to $\gamma_{0, C}$ is given by the \emph{Cameron--Martin formula} \citep[Corollary~2.4.3]{Bogachev1998}
\begin{equation}
	\label{eq:Cameron--Martin_formula}
	\gamma_{h, C}(\rd x) = \exp \left( \widehat{h}(x) - \frac{1}{2} \norm{h}_{\CMspace{\gamma_{0, C}}}^{2} \right) \, \gamma_{0, C}(\rd x),
\end{equation}
where $\widehat{h}(x)$ is the unique \emph{measurable linear functional} on $\cH$ (i.e., the unique map in the $L^{2}(\gamma_{0, C})$-closure of the topological dual $\cH^{\ast}$) extending the Cameron--Martin inner product $\innerprod{h}{\quark}_{\CMspace{\gamma_{0, C}}}$ \citep[Theorems~2.10.9 and 2.10.11]{Bogachev1998}.
We use \eqref{eq:Cameron--Martin_formula} to prove the following lower bound on the Gaussian mass of balls:\footnote{The proof extends immediately to the case that $\cH$ is a Banach space; see \citet[Proposition~6.29]{Lambley2026}.}

\begin{lemma}
	\label{lem:Gaussian_ball_mass_lower_bound}
	Let $\cH$ be a separable Hilbert space and let $\gamma_{0, C}$ be a centred, non-degenerate Gaussian measure on $\cH$.
	Then, for any $h \in \CMspace{\gamma_{0, C}}$ and any $r > 0$,
	\begin{equation*}
		\gamma_{0, C}\bigl( \oball{h}{r} \bigr) \geq \exp\left(-\frac{1}{2} \norm{h}_{\CMspace{\gamma_{0, C}}}^{2} \right) \gamma_{0, C}\bigl( \oball{0}{r} \bigr).
	\end{equation*}
\end{lemma}

\begin{proof}
	By the Cameron--Martin formula \eqref{eq:Cameron--Martin_formula},
	\begin{equation} \label{eq:Gaussian_lower_bound_CM}
		\gamma_{0, C}\bigl( \oball{h}{r} \bigr) = \gamma_{-h, C}\bigl(\oball{0}{r} \bigr) = \int_{\oball{0}{r}} \exp\left( -\frac{1}{2} \norm{h}^{2}_{\CMspace{\gamma_{0, C}}} - \widehat{h}(x) \right) \, \gamma_{0, C}(\rd x).
	\end{equation}
	Since $\gamma_{0, C}$ is symmetric, $\int_{\oball{0}{r}} \widehat{h}(x) \,\gamma_{0, C}(\rd x) = 0$ \citep[see][Lemma~4.7.1]{Bogachev1998}.
	Thus, by Jensen's inequality,
	\begin{align*}
		1 &= \exp\left(\frac{1}{\gamma_{0, C}(\oball{0}{r})} \int_{\oball{0}{r}} -\widehat{h}(x) \,\gamma_{0, C}(\rd x)\right) \\
		&\leq \frac{1}{\gamma_{0, C}(\oball{0}{r})} \int_{\oball{0}{r}} \exp\bigl( - \widehat{h}(x) \bigr) \,\gamma_{0, C}(\rd x).
	\end{align*}
	Applying the above in \eqref{eq:Gaussian_lower_bound_CM} yields the claim.
\end{proof}

\Cref{lem:Gaussian_ball_mass_lower_bound} is used only to prove the following result on the asymptotics of Gaussian ball masses, or equivalently of the function $F(r)$:

\begin{lemma}
	\label{lem:Gaussian_rapid_variation}
	Let $\cH$ be an infinite-dimensional separable Hilbert space and let $\gamma_{0, C}$ be a centred, non-degenerate Gaussian measure on $\cH$.
	Then, for any $s \in (0,1)$,
	\begin{equation} \label{eq:rapid_variation_claim}
		\lim_{r \to 0} \frac{\gamma_{0, C}(\oball{0}{s r})}{\gamma_{0, C}(\oball{0}{r})} = 0.
	\end{equation}
\end{lemma}

\begin{proof}
	Let $N \in \bN$ and let $(v_{n})_{n = 1}^{N} \subseteq \oball{0}{\tfrac{3}{4}}$ satisfy $\norm{v_{n} - v_{m}} > \tfrac{1}{2}$ for all $n \neq m$;
	such a sequence exists by Riesz's lemma as $\cH$ has infinite dimension.
	Since $\gamma_{0, C}$ is non-degenerate, $\CMspace{\gamma_{0, C}}$ is dense in $\cH$, and so we may arbitrarily select a sequence $(h_{n})_{n = 1}^{N} \subset \oball{0}{\tfrac{3}{4}} \cap \CMspace{\gamma_{0, C}}$ satisfying $\norm{h_{n} - v_{n}} < \tfrac{1}{16}$ for all $1 \leq n \leq N$ and $\norm{h_{n} - h_{m}} > \tfrac{1}{4}$ for all $n \neq m$.
	Then, by \Cref{lem:Gaussian_ball_mass_lower_bound}, for any $r \in (0, 1)$ we have
	\begin{equation*}
		\gamma_{0, C}\bigl(\oball{0}{r}\bigr)
		\geq \sum_{n = 1}^{N} \gamma_{0, C}\bigl(\Oball{r h_{n}}{\tfrac{r}{8}}\bigr)
		\geq \sum_{n = 1}^{N} \exp \left( - \frac{r^{2}}{2} \norm{h_{n}}_{\CMspace{\gamma_{0, C}}}^{2} \right) \gamma_{0, C}\bigl(\Oball{0}{\tfrac{r}{8}}\bigr).
	\end{equation*}
	From this, we obtain that
	\begin{equation*}
		0 \leq \limsup_{r \to 0} \frac{\gamma_{0, C}(\oball{0}{\frac{r}{8}})}{\gamma_{0, C}(\oball{0}{r})} \leq \limsup_{r \to 0} \left(\sum_{n = 1}^{N} \exp \left( - \frac{r^{2}}{2} \norm{h_{n}}_{\CMspace{\gamma_{0, C}}}^{2} \right) \right)^{-1} = N^{-1}.
	\end{equation*}
	Since the above holds for all $N \in \bN$, the claim \eqref{eq:rapid_variation_claim} holds for all $0 < s \leq \tfrac{1}{8}$.
	For $\tfrac{1}{8} < s < 1$, the logarithmic concavity of $\gamma_{0, C}$ \citep[Theorem~4.2.2]{Bogachev1998} implies that
	\begin{align*}
		\gamma_{0, C}\bigl( \oball{0}{r} \bigr)
		&= \gamma_{0, C}\bigl( \eta \oball{0}{s r} + (1 - \eta) \oball{0}{8s r} \bigr) \\
		&\geq \gamma_{0, C}\bigl( \oball{0}{s r} \bigr)^{\eta} \gamma_{0, C}\bigl(\oball{0}{8s r}\bigr)^{1-\eta},
	\end{align*}
	where $\eta = \tfrac{8}{7} - \tfrac{1}{7} s^{-1} > 0$.
	Thus,
	\begin{align*}
		\limsup_{r \to 0} \frac{\gamma_{0, C}(\oball{0}{s r})}{\gamma_{0, C}(\oball{0}{r})} &\leq \limsup_{r \to 0} \frac{\gamma_{0, C}(\oball{0}{s r})}{\gamma_{0, C}(\oball{0}{s r})^{\eta} \gamma_{0, C}(\oball{0}{8s r})^{1-\eta}}\\
		&= \limsup_{r \to 0} \left(\frac{\gamma_{0, C}(\oball{0}{s r})}{\gamma_{0, C}(\oball{0}{8s r})}\right)^{1-\eta} = 0. \qedhere
	\end{align*}
\end{proof}

\begin{remark}
	\begin{enumerate}[label=(\alph*)]
		\item
		The property \eqref{eq:rapid_variation_claim} asserts that ball probabilities of any non-degenerate Gaussian measure on $\cH$ are \emph{rapidly varying at $r = 0$} \citep[see][Section~2.4]{BinghamGoldieTeugels1987} when $\cH$ has infinite dimension.
		When $\cH = \bR^{d}$ has finite dimension, this is not the case:
		\begin{equation*}
			\lim_{r \to 0} \frac{\gamma_{0, C}\bigl( \oball{0}{sr} \bigr)}{\gamma_{0, C}\bigl(\oball{0}{r}\bigr)} = \lim_{r \to 0} \frac{\rho_{0,C}(0) \mathcal{L}^{d}\bigl(\oball{0}{sr} \bigr)}{\rho_{0,C}(0) \mathcal{L}^{d}\bigl( \oball{0}{r} \bigr)} = s^{d} > 0,
		\end{equation*}
		where $\mathcal{L}^{d}$ denotes $d$-dimensional Lebesgue measure on $\bR^{d}$.

		\item
		The rapid variation is related to the fact that no infinite-dimensional Hilbert space admits a probability measure $\mu$ that is \emph{doubling} \citep[Theorem~1]{Ruzhansky2001}, i.e., there can be no $\mu$ satisfying, for some $A > 0$,
		\begin{equation*}
			\frac{\mu\bigl(\oball{x}{\tfrac{r}{2}}\bigr)}{\mu\bigl(\oball{x}{r} \bigr)} \geq A > 0 \quad \text{for all $x \in \supp(\mu)$ and $r > 0$.}
		\end{equation*}
	\end{enumerate}
\end{remark}

As a final technical result, we bound the eigenvalues of the covariance operators $C_{v}$ defined in \eqref{eq:disintegration_covariance}. 
This captures the idea that the direction $\psi_{1}$ has the largest variance under the Gaussian measure $\gamma_{0,C}$;
thus, the conditional covariance $C_{\psi_{1}}$ must have smaller eigenvalues than any other such operator $C_{v}$ since the direction of greatest variance has been conditioned upon:

\begin{lemma}[Bounds on eigenvalues of conditioned Gaussian covariance]
	\label{lem:min-max}
	Let $\cH$ be an infinite-dimensional separable Hilbert space and let $C \colon \cH \to \cH$ be a trace-class \ac{SPSD} operator with eigenvalues $\lambda_{n}(C) \geq 0$ and corresponding unit eigenvectors $\psi_{n}$, $n \in \bN$.
	Let $v \in \cS \subset \cH$ have $\sigma_{v}^{2} \defeq \innerprod{ C v }{ v } > 0$, and let $C_{v} \colon \perpspace{\{v\}} \to \perpspace{\{v\}}$ be defined as in \eqref{eq:disintegration_covariance}.
	Then
	\begin{equation*}
		\lambda_{n}\bigl(C_{\psi_{1}}\bigr) \leq \lambda_{n}\bigl(C_{v}\bigr) \quad \text{for all $n \in \bN$.}
	\end{equation*}
\end{lemma}

\begin{proof}
	Recall that $C_{v} = P_{\perpspace{\{v\}}} R_{v} P^{\ast}_{\perpspace{\{v\}}}$ with $R_{v} \defeq C - \innerprod{Cv}{v}^{-1} (Cv \otimes Cv)$ as defined in \eqref{eq:disintegration_covariance_R}.
	Since $R_{v} = C_{v} \oplus 0$ on $\cH = \perpspace{\{v\}} \oplus \Span\{v\}$, it follows that $\lambda_{n}(C_{v}) = \lambda_{n}(R_{v})$ for all $n \in \bN$.
	Moreover, $\lambda_{n}(C_{\psi_{1}}) = \lambda_{n}(R_{\psi_{1}}) = \lambda_{n+1}(C)$ for all $n \in \bN$.

	By the min-max principle \eqref{eq:Courant--Fischer--Weyl_min-max} for the eigenvalues of the compact \ac{SPSD} operator $R_{v}$,
	\begin{equation}
		\label{eq:Courant--Fischer--Weyl_min-max_specialised}
		\lambda_{n}(R_{v}) = \min_{\substack{W_{n-1} \subseteq \cH\\\dim W_{n-1} = n-1}} \max_{x \in W_{n-1}^{\perp} \cap \cS} \innerprod{R_{v} x}{x}.
	\end{equation}
	For each $n \in \bN$, let $W_{n-1} \subseteq \cH$ be any $(n-1)$-dimensional subspace of $\cH$ attaining the minimum in \eqref{eq:Courant--Fischer--Weyl_min-max_specialised}, and let $V_{n}$ be any $n$-dimensional subspace of $\cH$ containing $W_{n-1} + \Span\{Cv\}$.
	Then, for any $x \in \perpspace{V_{n}} \subseteq \perpspace{\bigl(\Span\{Cv\}\bigr)}$, we see that
	\begin{equation} \label{eq:Rv_equals_Cv_on_perpspace}
		R_{v} x = Cx - \innerprod{Cv}{v}^{-1} \underset{=0}{\underbrace{\innerprod{Cv}{x}}} Cv = Cx.
	\end{equation}
	Applying the min-max principle on $C$, and using the identities \eqref{eq:Rv_equals_Cv_on_perpspace} and \eqref{eq:Courant--Fischer--Weyl_min-max_specialised}, we see that
	\begin{equation*}
		\lambda_{n}(C_{\psi_{1}}) = \lambda_{n+1}(C) 
		\leq  \max_{x \in \perpspace{V_{n}} \cap \cS} \innerprod{R_{v}x}{x} \leq \max_{x \in \perpspace{W_{n-1}} \cap \cS} \innerprod{R_{v}x}{x} = \lambda_{n}(R_{v}) = \lambda_{n}(C_{v}). \qedhere
	\end{equation*}
\end{proof}

\section{Discussion}
\label{sec:discussion}

In applications to \acp{BIP} and the analysis of random dynamical systems, one is often interested in computing the `most likely points' of a Borel probability measure $\mu^{\Phi}$ of the form
\begin{equation}
	\label{eq:posterior}
	\mu^{\Phi}(\rd x) = Z^{-1} \exp \bigl( - \Phi(x) \bigr) \, \mu(\rd x) , \qquad 0 < Z \defeq \int_{X} \exp\bigl(-\Phi(x)\bigr) \,\mu(\rd x) < \infty,
\end{equation}
where $\mu$ is a reference probability measure on a metric space $X$.
In the Bayesian context, $\mu$ is called the \emph{prior}, the function $\Phi$ is called the \emph{negative log-likelihood} or \emph{potential}, the measure $\mu^{\Phi}$ is called the \emph{posterior}, and its `most likely points' are called \ac{MAP} estimators.
In this setting, a notion of \ac{MAP} estimator is needed that does not rely on (Lebesgue) probability density functions:
such densities do not exist when $X$ is an infinite-dimensional Hilbert or Banach space, as is frequently the case in the nonparametric approach to \acp{BIP}.
Instead, a \ac{MAP} estimator is understood as the centre $u$ of a ball of asymptotically maximal posterior mass in the limit as ball radius tends to zero, as can be made precise in the following way:

\begin{definition}
	Let $X$ be a separable metric space and let $\mu$ be a Borel probability measure on $X$.
	Then $x \in X$ is a \emph{strong mode} for $\mu$ if
	\begin{equation*}
		\lim_{r \to 0} \frac{\mu\bigl(B(x, r)\bigr)}{\sup_{y \in X} \mu\bigl(B(y, r)\bigr)} = 1.
	\end{equation*}
	Similarly, $x \in X$ is a \emph{weak mode} for $\mu$ if, for all $y \in \supp(\mu)$,
	\begin{equation*}
		\liminf_{r \to 0} \frac{\mu\bigl( B(x, r) \bigr)}{\mu\bigl( B(y, r) \bigr)} \geq 1.
	\end{equation*}
	Strong and weak modes of a measure $\mu^{\Phi}$ of the form \eqref{eq:posterior} are often called \emph{strong \ac{MAP} estimators} and \emph{weak \ac{MAP} estimators} respectively.
\end{definition}

Any strong mode is clearly also a weak mode, but the converse need not be true.
Weak modes are closely connected to \ac{OM} functionals:

\begin{proposition}[{\citealp[Proposition~4.1]{AyanbayevKlebanovLieSullivan2022I}}]
	\label{prop:OM_weak}
	Let $X$ be a separable metric space and let $\mu$ be a Borel probability measure on $X$.
	Suppose that $\mu$ admits  an exhaustive \ac{OM} functional $I_{\mu} \colon X \to \bR \cup \{+\infty\}$.
	Then the weak modes of $\mu$ are precisely the minimisers of $I_{\mu}$.
\end{proposition}

In combination with the following \emph{transfer lemma}, \Cref{prop:OM_weak} gives a straightforward method for proving the existence of weak \ac{MAP} estimators for measures of the form \eqref{eq:posterior}:

\begin{lemma}[{Transfer lemma for exhaustive \ac{OM} functionals; \citealp[Proposition~2.6]{Lambley2023}}]
	\label{lem:transfer_lemma}
	Let $X$ be a separable metric space, let $\mu$ be a Borel probability measure on $X$, and let $\mu^{\Phi}$ be given by \eqref{eq:posterior} for some continuous potential $\Phi \colon X \to \bR$.
	If $\mu$ has an exhaustive \ac{OM} functional $I_{\mu} \colon X \to \bR \cup \{+\infty\}$, then $\mu^{\Phi}$ has an exhaustive \ac{OM} functional $I_{\mu^{\Phi}} \colon X \to \bR \cup \{+\infty\}$ given by
	\begin{equation*}
		I_{\mu^{\Phi}}(x) \defeq I_{\mu}(x) + \Phi(x).
	\end{equation*}
\end{lemma}

Given that, by \Cref{ex:Gaussian_OM}, a Gaussian measure $\gamma_{0, C}$ on a separable Hilbert space $\cH$ has an exhaustive \ac{OM} functional $I_{\gamma_{0, C}}(x) = \tfrac{1}{2} \norm{x}_{\CMspace{\gamma_{0, C}}}^{2}$, \Cref{lem:transfer_lemma} immediately yields an exhaustive \ac{OM} functional for Gaussian measures reweighted by a continuous density as in  \eqref{eq:posterior}.
Moreover, the minimisers of this \ac{OM} functional are the weak \ac{MAP} estimators.

Deriving the strong \ac{MAP} estimators of such reweightings requires more work, and this has been studied  extensively in Hilbert spaces \citep{DashtiLawStuartVoss2013,Kretschmann2023} and in sequence spaces \citep{KlebanovWacker2023}, with the general Banach case addressed by \citet{Lambley2023}.
Indeed, under a mild lower-bound assumption on $\Phi$, the minimisers of this posterior \ac{OM} functional are also strong \ac{MAP} estimators \citep[Theorem~1.1]{Lambley2023}.
Thus, \ac{MAP} estimation with the Gaussian prior $\gamma_{0, C}$ coincides with Tikhonov regularisation over the Cameron--Martin space, with the Cameron--Martin norm providing the Tikhonov penalty:
\begin{align*}
	u \text{ is a \ac{MAP} estimator of } \gamma_{0, C}^{\Phi}
	&\iff
	u \text{ minimises } \Phi + \frac{1}{2} \norm{ \quark }_{\CMspace{\gamma_{0, C}}}^{2} \text{ in } \cH \\
	&\iff
	u \text{ minimises } \Phi + \frac{1}{2} \norm{ \quark }_{\CMspace{\gamma_{0, C}}}^{2} \text{ in } \CMspace{\gamma_{0, C}} .
\end{align*}

On the basis of the Gaussian case, one might expect that the computation of \ac{MAP} estimators under an \ac{ACG} prior $\nu_{C}$ involves the minimisation of $\Phi$, or a penalised version of $\Phi$, over $\cS$ or perhaps $\CMspace{\gamma_{0, C}} \cap \cS$.
However, the results of this paper show that \ac{MAP} estimators under an \ac{ACG} prior are in fact the minimisers of a much more tightly constrained problem:

\begin{corollary}
	\label{cor:OM_of_reweighted_ACG}
	Let $\nu_{C}$ be a non-degenerate \ac{ACG} measure on the unit sphere $\cS$ of an infinite-dimensional, separable Hilbert space $\cH$.
	Let $L = \ker(C - \lambda_{1}(C) I)$ be the subspace of $\cH$ spanned by the eigenvectors of $C$ with maximal eigenvalue.
	Let $\Phi \colon \cS \to \bR$ be continuous and suppose that $Z \defeq \int_{\cS} e^{-\Phi(x)} \, \nu_{C}(\rd x) < \infty$.
	Let $\nu_{C}^{\Phi}$ be the probability measure on $\cS$ given by
	\[
		\nu_{C}^{\Phi}(\rd x) = Z^{-1} \exp \bigl( - \Phi(x) \bigr) \, \nu_{C}(\rd x) .
	\]
	\begin{enumerate}[label=(\alph*)]
		\item
		\label{item:OM_posterior}
		The measure $\nu_{C}^{\Phi}$ has an exhaustive \ac{OM} functional $I_{\nu_{C}^{\Phi}} \colon \cS \to \bR \cup \{ + \infty \}$ given by
		\begin{equation}
			\label{eq:reweighted_OM_ACG}
			I_{\nu_{C}^{\Phi}}(x) \defeq I_{\nu_{C}}(x) + \Phi(x) =
			\begin{cases}
				\Phi(x), & \text{if $x \in \cS \cap L$,} \\
				+\infty, & \text{if $x \in \cS \setminus L$.}
			\end{cases}
		\end{equation}

		\item 
			\label{item:weak_modes_and_posterior_OM_minimisers}
		The set of weak \ac{MAP} estimators of $\nu_{C}^{\Phi}$ is non-empty and coincides with the set of minimisers of $I_{\nu_{C}^{\Phi}}$.

		\item
		\label{item:weak_and_strong_modes}
		If $\Phi$ is bounded below, then the sets of strong and weak \ac{MAP} estimators of $\nu_{C}^{\Phi}$ coincide.
	\end{enumerate}
\end{corollary}

\begin{proof}
	By \Cref{thm:main}, the function $I_{\nu_{C}}$ defined in \eqref{eq:main} is an exhaustive \ac{OM} functional for $\nu_{C}$.
	Then, by \Cref{lem:transfer_lemma}, $I_{\nu_{C}^{\Phi}}$ is given by the sum of $I_{\nu_{C}}$ and $\Phi$, yielding \eqref{eq:reweighted_OM_ACG}.
	For \ref{item:weak_modes_and_posterior_OM_minimisers}, the fact that the weak \ac{MAP} estimators coincide with the minimisers of $I_{\nu_{C}^{\Phi}}$ follows from \Cref{prop:OM_weak}, and the minimising set is non-empty as $\Phi$ is continuous and $\cS \cap L$ is compact.

	For \ref{item:weak_and_strong_modes}, we show that all weak \ac{MAP} estimators are strong \ac{MAP} estimators.
	Since all strong \ac{MAP} estimators are weak \ac{MAP} estimators, this will show that the strong and weak \ac{MAP} estimators coincide.
	Thus, let $v^{\star} \in \argmin_{v \in \cS \cap L} \Phi(v)$ be any weak \ac{MAP} estimator and fix positive null sequences $(\theta_{j})_{j \in \bN}$ and $(\varepsilon_{j})_{j \in \bN}$.
	Suppose that $(w_{j})_{j \in \bN} \subset \cS$ is any sequence satisfying
	\begin{equation}
		\label{eq:asymptotically_maximising}
		\nu_{C}^{\Phi}\bigl( \angball{w_{j}}{\theta_{j}} \bigr) \geq (1 - \varepsilon_{j}) \sup_{v \in \cS} \nu_{C}^{\Phi}\bigl( \angball{v}{\theta_{j}} \bigr).
	\end{equation}
	Then, using that $\Phi$ is bounded below and continuous in the first line, and applying \eqref{eq:asymptotically_maximising} in the second line, we have
	\begin{align*}
		\liminf_{j \to \infty} \frac{\nu_{C}\bigl(\angball{w_{j}}{\theta_{j}} \bigr)}{\nu_{C}\bigl(\angball{\psi_{1}}{\theta_{j}} \bigr)} &\geq \liminf_{j \to \infty} \frac{e^{\inf_{v} \Phi(v)} \nu_{C}^{\Phi}\bigl( \angball{w_{j}}{\theta_{j}} \bigr)}{e^{\Phi(\psi_{1})} \nu_{C}^{\Phi} \bigl( \angball{\psi_{1}}{\theta_{j}} \bigr)} \\
		&\geq \exp\Bigl( \inf_{v \in \cS} \Phi(v) - \Phi(\psi_{1}) \Bigr) \liminf_{j \to \infty} \frac{\nu_{C}^{\Phi}\bigl( \angball{w_{j}}{\theta_{j}} \bigr)}{\sup_{v \in \cS} \nu_{C}^{\Phi}\bigl( \angball{v}{\theta_{j}} \bigr)} > 0.
	\end{align*}
	Thus, by \eqref{eq:acg_exhaustivity}, we must have $\sigma_{w_{j}}^{2} \defeq \innerprod{Cw_{j}}{w_{j}} \to \innerprod{C\psi_{1}}{\psi_{1}} \qefed \sigma_{\psi_{1}}^{2}$.
	This implies that $(w_{j})_{j \in \bN}$ has a convergent subsequence.
	To see this, write $w_{j} = P_{L} w_{j} + P_{L^{\perp}} w_{j}$ and note that, since $L$ has finite dimension and $\norm{P_{L} w_{j}} \leq \norm{w_{j}} \leq 1$, the sequence $(P_{L} w_{j})_{j \in \bN}$ must converge along some strictly increasing sequence $(j_{k})_{k \in \bN} \to \infty$;
	call the limit $w^{\star}$.
	Next, using that $\norm{P_{L}w_{j}}^{2} = \norm{w_{j}}^{2} - \norm{P_{L^{\perp}} w_{j}}^{2} = 1 - \norm{P_{L^{\perp}} w_{j}}^{2}$, as well as that $\norm{C P_{L}} = \lambda_{1}(C)$ and $\norm{C P_{L^{\perp}}} < \lambda_{1}(C)$, 
	\begin{align*}
		\innerprod{Cw_{j}}{w_{j}} = \innerprod{C P_{L} w_{j}}{P_{L} w_{j}} + \innerprod{C P_{L^{\perp}} w_{j}}{P_{L^{\perp}}w_{j}} &\leq \norm{C P_{L}} \norm{P_{L} w_{j}}^{2} + \norm{C P_{L^{\perp}}} \norm{P_{L^{\perp}}w_{j}}^{2} \\
		&= \lambda_{1}(C) - \underbrace{(\lambda_{1}(C) - \norm{C P_{L^{\perp}}} )}_{> 0} \norm{P_{L^{\perp}} w_{j}}^{2}.
	\end{align*}
	Taking $j \to \infty$, we conclude that $\norm{P_{L^{\perp}} w_{j}} \to 0$, and so $(w_{j_{k}})_{k \in \bN} \to w^{\star}$.
	Then, since $w^{\star} \in \cS \cap L$, we have by \eqref{eq:mass_same_in_leading_eigenspace} and \eqref{eq:acg_exhaustivity_bound_2} that 
	\begin{align}
		\nu_{C}(\angball{w^{\star}}{\theta}) &= \sup_{v \in \cS} \nu_{C}(\angball{v}{\theta}) &&\text{for all $0 < \theta < \tfrac{\pi}{2}$.} \label{eq:w_supremal}
	\end{align}
	Thus, we conclude that
	\begin{align*}
		& \liminf_{k \to \infty} \frac{\nu_{C}^{\Phi}\bigl( \angball{v^{\star}}{\theta_{j_{k}}} \bigr)}{\sup_{v \in \cS} \nu_{C}^{\Phi}\bigl( \angball{v}{\theta_{j_{k}}} \bigr)} \\
		& \quad \geq \underset{\geq 1\text{~using $I_{\nu_{C}^{\Phi}}$}}{\underbrace{\liminf_{k \to \infty} \frac{\nu_{C}^{\Phi}\bigl( \angball{v^{\star}}{\theta_{j_{k}}} \bigr)}{\nu_{C}^{\Phi}\bigl( \angball{w^{\star}}{\theta_{j_{k}}} \bigr)  } }}	\liminf_{k \to \infty} \frac{\nu_{C}^{\Phi}\bigl( \angball{w^{\star}}{\theta_{j_{k}}} \bigr)}{\sup_{v \in \cS} \nu_{C}^{\Phi}\bigl( \angball{v}{\theta_{j_{k}}} \bigr)} \\
		& \quad \geq \liminf_{k \to \infty} \frac{\nu_{C}^{\Phi}\bigl(\angball{w^{\star}}{\theta_{j_{k}}}\bigr)}{\nu_{C}^{\Phi}\bigl(\angball{w_{j_{k}}}{\theta_{j_{k}}} \bigr)}  \underset{=1 \text{~by \eqref{eq:asymptotically_maximising}}}{\underbrace{\liminf_{k \to \infty} \frac{\nu_{C}^{\Phi}\bigl( \angball{w_{j_{k}}}{\theta_{j_{k}}} \bigr)}{\sup_{v \in \cS}\nu_{C}^{\Phi}\bigl( \angball{v}{\theta_{j_{k}}}\bigr)}}}\\
		& \quad \geq \liminf_{k \to \infty} \frac{Z^{-1} \inf_{\angball{w^{\star}}{\theta_{j_{k}}}} e^{-\Phi} \nu_{C}\bigl(\angball{w^{\star}}{\theta_{j_{k}}} \bigr)}{Z^{-1} \sup_{\angball{w_{j_{k}}}{\theta_{j_{k}}}} e^{-\Phi} \nu_{C}\bigl( \angball{w_{j_{k}}}{\theta_{j_{k}}}\bigr)} &&\text{(using definition of $\nu_{C}^{\Phi}$)} \\
		& \quad \geq \frac{e^{-\Phi(w^{\star})}}{e^{-\Phi(w^{\star})}} =  1. &&\text{(by \eqref{eq:w_supremal} and continuity of $\Phi$)}
	\end{align*}
	Since, for every null sequence $(\theta_{j})_{j \in \bN}$ of radii, the claim holds along a subsequence, it follows from Urysohn's subsequence principle that
	\begin{equation*}
		\lim_{\theta \to 0} \frac{\nu_{C}^{\Phi}\bigl(\angball{v^{\star}}{\theta} \bigr)}{\sup_{v \in \cS} \nu_{C}^{\Phi} \bigl( \angball{v}{\theta} \bigr)} = 1.
	\end{equation*}
	This shows that any weak \ac{MAP} estimator $v^{\star}$ is also strong, completing the proof.
\end{proof}

By the strong--weak dichotomy \citep[Lemma~2.2]{Lambley2023}, either the strong and weak \ac{MAP} estimators coincide, or the set of strong \ac{MAP} estimators is empty.
As discussed earlier, the use of a lower bound on $\Phi$ to ensure that a strong \ac{MAP} estimator exists is standard for Gaussian priors, so it is not surprising that such a bound is also useful for \ac{ACG} priors.
We leave the question of whether strong \ac{MAP} estimators exist for potentials that are unbounded below to future work.

\Cref{cor:OM_of_reweighted_ACG} is somewhat troubling for Bayesian inverse problems with \ac{ACG} priors, since it shows that the computation of a \ac{MAP} estimator for such a problem is not the minimisation of the negative log-likelihood $\Phi$ with an \emph{additive penalty} coming from the negative log-prior, but is instead a \emph{constrained minimisation} of $\Phi$.
Moreover, this constrained minimisation is not restricted to the relatively rich set $\cS \cap \CMspace{\gamma_{0, C}}$ of unit vectors with finite Cameron--Martin norm, but to the much smaller set $\cS \cap L$ of leading eigenvectors of $C$.
Indeed, if the leading eigenvalue of $C$ is simple, then $\cS \cap L$ is a two-point set!
In addition, the \ac{MAP} estimators of such problems are completely insensitive to features of the prior below the level of the leading eigenvalue:
taking a different covariance operator where one has perturbed the remaining eigenvalues to any values in $(0, \lambda_{1}(C))$, or where one has altered the corresponding eigenvectors, does not change $\cS \cap L$ or the minimisers of $\Phi$ in $\cS \cap L$.

\section{Future work}
\label{sec:future}

In light of the somewhat degenerate behaviour of the \ac{ACG} measure in infinite dimensions, it is natural to ask whether other spherical distributions  \citep[see][p.~182]{MardiaJupp1999} are better behaved.
It would be interesting to know, for example, whether the projections of other measures onto the infinite-dimensional sphere are similarly degenerate.
To adapt the proof strategy used here, one would require an analogue of the rapid-variation property (\Cref{lem:Gaussian_rapid_variation}) and the conditioning formula for Gaussian measures in infinite dimensions.

\section*{Acknowledgements and declarations}
\addcontentsline{toc}{section}{Acknowledgements}

HL was supported by the Warwick Mathematics Institute Centre for Doctoral Training and gratefully acknowledges funding from the University of Warwick and the UK Engineering and Physical Sciences Research Council (Grant number: EP/W524645/1).
HL is also funded by the Deutsche Forschungsgemeinschaft (DFG, German Research Foundation) under Germany's Excellence Strategy --- The Berlin Mathematics Research Center MATH+ (EXC-2046/2, project ID: 390685689).

For the purpose of open access, the authors have applied a Creative Commons Attribution (CC~BY) licence to any Author Accepted Manuscript version arising.
No data relating to humans, living or dead, were created or analysed in this study.

The manuscript was written by the named authors.
ChatGPT~6 Astra, Claude Sonnet~5.5, and Gemini~3.1 Pro were used before submission to identify typographical errors and mathematical omissions, which were then corrected manually.
ChatGPT~5.5 was used to investigate general decompositions for the integral in \eqref{eq:nu_upper_bound}, which led the authors to adopt the approach of \Cref{lem:Gaussian_rapid_variation}.

\bibliographystyle{abbrvnat}
\bibliography{references}
\addcontentsline{toc}{section}{References}

\clearpage
\appendix

\end{document}